\documentclass[12pt]{amsart}
\usepackage{rotating,pdflscape,color,amssymb,hyperref,subfigure,psfrag,amsmath,bbm,epsfig}
\usepackage{a4wide,amscd}
\usepackage{tikz} 
 \usepackage[usenames,dvipsnames]{pstricks}
 \usepackage{pst-grad} 
 \usepackage{pst-plot} 

\usepackage{amsmath}
\usepackage{amsfonts}
\usepackage{amssymb}

\usepackage[small]{caption}

\newcommand{\ld}{\ldots}
\newcommand{\beg}{\begin}
\newcommand{\en}{\end}

\newcommand{\cdd}{\cdots\cdots}

\newcommand{\trm}{\textrm}

\newcommand{\bgt}{\begin{itemize}}
\newcommand{\ent}{\end{itemize}}
\newcommand{\ite}{\item}

\newcommand{\op}{\operatorname}
\newcommand{\eqre}{\eqref}
\newcommand{\re}{\ref}
\newcommand{\la}{\label}

\newcommand{\rfl}{\rfloor}
\newcommand{\lfl}{\lfloor}

\newcommand{\h}{\mc{H}}

\newcommand{\diag}{\operatorname{diag}}

\newcommand{\p}{\mathbb{P}}

\newcommand{\Tr}{\operatorname{Tr}}

\newcommand{\ssi}{if and only if }

\newcommand{\z}{\mathbb{Z}}

\newcommand{\ud}{\mathrm{d}}

\newcommand{\f}{\frac}
\newcommand{\ff}{\frac{1}}
\newcommand{\lf}{\left}
\newcommand{\ri}{\right}

\newcommand{\st}{such that }

\newcommand{\lam}{\lambda}

\newcommand{\ti}{\times}

\newcommand{\ste}{\, ;\, }
\newcommand{\mc}{\mathcal }

\newcommand{\A}{\mc{A}}

\newcommand{\al}{\alpha}

\newcommand{\eqlaw}{\stackrel{\textrm{law}}{=}}

\newcommand{\bbm}{\begin{bmatrix}}
\newcommand{\ebm}{\end{bmatrix}}
\newcommand{\bes}{\begin{equation*}}
\newcommand{\ees}{\end{equation*}}
\newcommand{\be}{\begin{equation}}
\newcommand{\ee}{\end{equation}}
\newcommand{\beqy}{\begin{eqnarray}}
\newcommand{\eeqy}{\end{eqnarray}}
\newcommand{\beq}{\begin{eqnarray*}}
\newcommand{\eeq}{\end{eqnarray*}}
\newcommand{\one}{\mathbbm{1}}

\newcommand{\ie}{\emph{i.e. }}

\newcommand{\bpm}{\begin{pmatrix}}
\newcommand{\epm}{\end{pmatrix}}

\newcommand{\cd}{\cdots}

\newcommand{\bpr}{\beg{proof}}
\newcommand{\epr}{\en{proof}}

\newcommand{\del}{\delta}
\newcommand{\Del}{\Delta}

\newcommand{\Pa}{\mc{P}}

\newcommand{\wtH}{\widetilde{H}} 
\newcommand{\tra}{^\top}

\theoremstyle{definition}

\long\def\symbolfootnote[#1]#2{\begingroup
\def\thefootnote{\fnsymbol{footnote}}\footnote[#1]{#2}\endgroup}

\begin{document}

\title[]{GUE minors,  
maximal Brownian functionals and longest increasing subsequences}
\author{Florent Benaych-Georges
\and
Christian Houdr\'e}\thanks{FBG: MAP5,
Universit\'e Paris Descartes,
45, rue des Saints-P\`eres
75270 Paris Cedex 06, France. florent.benaych-georges@parisdescartes.fr.
This work was done in part during this author's stay  at Georgia Institute of Technology,
Atlanta.\\
CH: {School of Mathematics, Georgia Institute of Technology,
Atlanta, Georgia, 30332-0160, USA.
houdre@math.gatech.edu.} This work was partially supported by the grant \#246283 from
the Simons Foundation and by a Simons Foundation Fellowship, grant \#267336.  Many thanks to the 
Laboratory 
MAS of \'Ecole Centrale, Paris, and to the LPMA of the 
Universit\'e Pierre et Marie Curie for their hospitality while 
part of this research was carried out.}



\beg{abstract} We present equalities in law between the spectra of the minors of a GUE matrix 
and some maximal functionals of independent Brownian motions.  In turn, these results 
allow to recover the limiting shape (properly centered and scaled) of the RSK Young diagrams 
associated with a random word as a function of the spectra of these minors. 
Since the length of the top row of the diagrams is the length of the longest increasing 
subsequence of the random word, the corresponding limiting law also follows. 
\en{abstract}

\maketitle

\section{Introduction}
It is by now well known that there exist strong and interesting connections between
directed percolation and random matrices.  The precise results we have in 
mind have their origins in 
the identity in law, due to Baryshnikov~\cite{baryshnikov} and Gravner, 
Tracy and Widom \cite{GTW}, between 
the maximal eigenvalue of an $M\ti M$ element of the GUE and a certain maximal functional of 
standard $M$-dimensional 
Brownian motion originating in queuing theory, with Glynn and Whitt~\cite{GW}.  
This first result has seen many extensions and 
complements.   For example, O'Connell and Yor~\cite{OY} as well as 
Bougerol and Jeulin~\cite{BJ} obtained identities in law between (different) multivariate 
Brownian functionals and the spectrum of the GUE whose equivalence is 
shown in Biane, Bougerol and  O'Connell~\cite{BBO}.  
Various related representations have also been 
put forward and studied for instance in Doumerc~\cite{doumerc}, 
Johansson~\cite{J1,J2}, O'Connell~\cite{O} 
to name but a few authors and pieces of work.

Our interest in such representations comes from the identification 
by Its, Tracy and Widom~\cite{TW, ITW1, ITW2} 
of the limiting length (properly centered and scaled) 
of the longest increasing subsequence of a random word whose size  tends to infinity as the maximal 
eigenvalue of a certain random matrix.   For example, 
in the case of a word with i.i.d.~uniformly distributed letters in  an alphabet of size 
$M$, the limiting law is the maximal eigenvalue of the $M\ti M$ traceless GUE.  
Moreover, positively answering a conjecture of Tracy and Widom, 
Johansson~\cite{KJAOM2001} showed that   the whole normalized limiting shape of 
the RSK Young diagrams associated with the random word is the spectrum 
of the $M\ti M$ traceless GUE.  
Since the length of the top row of the diagrams is the length of the longest increasing 
subsequence of the random word the maximal eigenvalue result of \cite{TW} followed.

Limiting laws expressed in terms of maximal Brownian 
functionals are also obtained in \cite{HL1}.  
These last representations involve 
{\it dependent} Brownian motions and do not clearly recover the 
results of \cite{TW} or \cite{ITW1, ITW2}, which themselves are mainly derived  
by analytical techniques.  To resolve this issue, we provide below an 
extension of Baryshnikov's  result \cite{baryshnikov} on the identification of the 
multivariate law of the maximal eigenvalues of the principal minors of a GUE matrix 
with some maximal functionals of a standard multidimensional Brownian motion.  
This allows us to circumvent the analytical approach and provides a mixed 
combinatorial/probabilistic methodology to the solutions of these finite 
alphabet longest increasing subsequence problems.  
Our hope is that Theorem~\re{208121}, below, will also be helpful to  
fully identify eigenvalues of random matrices as the limiting laws in the corresponding Markov 
random word problems (see Kuperberg's Conjecture 7 in  \cite{kuperberg2002}).  
In the Markovian setting, the analytical 
methodology is lacking, in contrast to the probabilistic 
one, and to date the limiting laws are mainly only expressed as Brownian functionals.  
Indeed, the multivariate functional appearing in Theorem~\re{208121} is exactly the one giving 
the limit law of the shape of the RSK image of a Markov random word 
in \cite{HL2}, the only difference being 
that the Brownian motions in \cite{HL2} are correlated.  
This correlation issue in maximal functionals often amounts to adding a condition on the 
trace of the random matrix (as in \cite{TW, KJAOM2001, ITW1, ITW2}).  
However, for general Markov random words the full 
identification of these functionals via random matrices remains open.  
For Markov random words with cyclic and symmetric transition matrix, the longest 
increasing subsequence will be asymptotically identified to the eigenvalues of some random matrices  once we will have a more general version of Theorem \re{208121} below where the Brownian motions are correlated.
Our intuition is that to get such a generalization of  Theorem \re{208121}, one needs to 
consider the minors of more general random matrices, namely  Gaussian Hermitian matrices 
with general Gaussian vector as diagonal.
Besides providing the final touch to an essentially 
probabilistic proof of the random word asymptotics problem, 
our results also allow us to shed new lights on the queuing 
framework by providing, for example, 
joint limiting laws involving departing times and service times 
of individual customers.\\ \\  \\

\section{Statements and proofs of the results}

Throughout, fix a positive integer $M$ and 
consider: \bgt\ite an $M\ti M$ GUE matrix $H=[h_{ij}]$, \ie, a standard Gaussian 
variable on the space of  $M\ti M$ Hermitian matrices endowed with the  
Euclidean scalar product given 
by $X\cdot Y=\Tr XY$,\\
\ite an $M$-dimensional standard Brownian 
motion  $B=(B_k(t))_{t\in [0,1],k=1, \ld, M}$.
\ent

For each $k=1, \ld, M$, denote by \be\la{16812.20h07}\mu_{1}^k\ge \cd\ge \mu_k^k,\ee the 
eigenvalues of the principal $k\ti k$ minor of $H$.  Next, introduce the 
set $$\mc{P}:=\{\pi: [0,1]\to\{1, \ld, M\} \trm{ c\`adl\`ag, non-decreasing}\},$$
and for 
$\pi\in \Pa$, 
let $$\Del_{\pi} (B):=\int_0^1\ud B_{\pi(t)}(t)=\sum_{j=1}^M(B_j(t_j)-B_j(t_{j-1})),$$
where 
$0=t_0\le t_1\le \cd\le t_M=1$ 
are such that  
$$\pi(\cdot)=\sum_{j=1}^{M-1}j\ti\one_{[t_{j-1}, t_j)}(\cdot)+M\ti\one_{[t_{M-1},t_M]}(\cdot).$$
To complete our notations, for $\pi_1, \pi_2\in \Pa$, we write $\pi_1<\pi_2$ whenever $\pi_1(t)<\pi_2(t)$,
for all $t\in [0,1]$.   Let us now state our first result which, in particular, when $\ell = 1$ below,  
identifies the joint law of the maximal eigenvalue of the principal minors of $H$ and  
therefore  extends Theorem 0.7 of \cite{baryshnikov}.

\beg{Th}\la{208121} The following equality in law holds true:
$$
\lf(\sum_{i=1}^\ell \mu_i^k\ri)_{1\le \ell\le k\le M}\!\!\!\!\eqlaw\quad\lf( \sup\lf\{\sum_{i=1}^\ell
\Del_{\pi_i} (B)\ste \pi_1, \ld, \pi_\ell \in \Pa,  \pi_1< \cd< \pi_\ell\le k \ri\}
\ri)_{1\le \ell\le k\le M}$$
\en{Th}

This theorem  also has a process version (where the matrix $H$ is replaced by an Hermitian Brownian motion,  and $B$ is taken up to time $t$ and not to time $1$.

 In in the following proof, and throughout, 
$\Longrightarrow$ indicates convergence in distribution.

\bpr Let $w_{N,M}:=[w_{ij}]$ ($ i=1, \ld, N$,  $j=1, \ld, M$) be an array of i.i.d.~geometric random 
variables with parameter $q\in (0,1)$, i.e., with law $\sum_{k\ge 0} q^k(1-q)\del_k$.  
Such variables have mean $e:=q/(1-q)$ and variance $v:=q/(1-q)^2$.

Applying the RSK correspondence to $w_{N,M}$ (see e.g., \cite{moerbeke} for an introduction to the RSK  correspondence applied to arrays of integers) gives a pair $(P,Q)$ of Young diagrams with the same shape. 
Let us denote the shape of these Young diagrams by $$\lam_1^M\ge \cd\ge \lam_M^M.$$
The exponent $M$ is here to emphasize the dependence on the dimension $M$ of the GUE matrix $H$ (the dependence on $N$ is implicit). 
In the same way, one can of course define, for each $k=1,\ld, M$, the shape 
$$\lam_1^k\ge \cd\ge \lam_k^k,$$ 
of the Young diagrams obtained by applying 
the RSK correspondence to the array  $w_{N,k}$, which is the array $w_{N,M}$ where all but the first $k$   
columns  have been removed.  
Note that  \be\la{255141}\bpm\beg{array}{ccccccc}&&&\lam_1^1&&&\\
&&\lam^2_1&&\lam_2^2&&\\
&\lam^3_1&&\lam^3_2&&\lam_3^3&\\
&\ld&&\ld&&\ld&\\
\lam_1^M&&\ld&&\ld&&\lam_M^M
\en{array}\epm,\ee
is a Gelfand-Tsetlin pattern, i.e., satisfies the interlacing 
inequalities $\lam_i^k\ge \lam_i^{k-1}\ge \lam_{i+1}^k$ ($1\le i< k\le M$). This can be seen by the fact that if $(P_k,Q_k)$ denotes the pair  associated to  $w_{N,k}$ by the RSK correspondance, then $P_{k-1}$ can be deduced from $P_k$ by removing all the boxes filled with the number $k$.

Let us now define the random 
variables \be\la{1481218h04}\xi_i^k:=\f{\lam_i^k-eN}{\sqrt{vN}}.\ee  
Then  we have the following lemma.  \beg{lem}\la{1481218h07ter} As $N$ tends to 
infinity, $$ \bpm\beg{array}{ccccccc}&&&\xi_1^1&&&\\
&&\xi^2_1&&\xi_2^2&&\\
&\xi^3_1&&\xi^3_2&&\xi_3^3&\\
&\ld&&\ld&&\ld&\\
\xi_1^M&&\ld&&\ld&&\xi_M^M
\en{array}\epm\;\;\Longrightarrow\;\; \bpm\beg{array}{ccccccc}&&&\mu_1^1&&&\\
&&\mu^2_1&&\mu_2^2&&\\
&\mu^3_1&&\mu^3_2&&\mu_3^3&\\
&\ld&&\ld&&\ld&\\
\mu_1^M&&\ld&&\ld&&\mu_M^M
\en{array}\epm,$$
where the $\mu^k_i$s are the ones 
introduced in \eqre{16812.20h07}.
\en{lem}
\bpr This lemma is stated as Proposition~3.12 in \cite{NJEJP06} (with 
slightly different notation).  However, for the convenience of the reader, we shall say a few words about its proof:

By a well known property of the geometric law, conditionally to $S:=\sum_{i,j} w_{i,j}$, the law of $w_{N,M}$ is the uniform law on the set of $N\ti M$ integer arrays  summing up to $S$. But the RSK correspondence $w\mapsto (P(w),Q(w))$ is a bijection when defined on this set, and the information carried by the Gelfand-Tsetlin pattern  $\lam$ of \eqre{255141} is the same as the one carried by $P(w)$. It follows that  conditionally to the event $S=s$, the law of the Gelfand-Tsetlin pattern $\lam$  of \eqre{255141}  can be recovered as follows: first choose $\lam^M_1\ge\cd\ge  \lam_M^M$ according to its particular conditional law, and then $\lam$ is uniformly distributed on the set of   integer valued Gelfand-Tsetlin patterns with first row $\lam^M_1\ge\cd\ge  \lam_M^M$. As a consequence, by conditioning, it is clear that the unconditional law of $\lam$ can be recovered exactly 
in the same way.

Moreover, it can also be proved \cite[Propo. 4.7]{baryshnikov} that conditionally on 
the value of its first row $\mu_1^M\ge \cd\ge \mu_M^M$, the law 
of $(\mu_i^k)_{1\le i\le k\le M}$ is the uniform law of real valued 
Gelfand-Tsetlin patterns with first row $\mu_1^M\ge \cd\ge \mu_M^M$.

Thus as $M$ is fixed, up to standard analysis (see e.g. \cite[Lem. 2.5]{baryshnikov}), it suffices to prove that the convergence in distribution stated in the lemma is true when restricted to the first rows: 
\be\la{2551418h}(\xi_1^M,\ld\ld,\xi_M^M)\;\;\Longrightarrow\;\;(\mu_1^M,\ld\ld,\mu_M^M)\ee

To prove \eqre{2551418h}, one only needs to notice that the joint distribution of its right hand side 
has density  \be\la{2551419h}\ff{\trm{normalizing constant}} \prod_{1\le i<j\le M}(\mu_j^M-\mu_i^M)^2e^{-\ff{2}\sum_i (\mu_i^M)^2}\ee (this is a well known fact, see e.g.  \cite[Th. 2.5.2]{alice-greg-ofer}), whereas it has been proved   by Johansson in \cite{J1}  that the joint density of the left hand side of \eqre{2551418h} with respect to the counting measure is \be\la{2551419h1} \ff{\trm{normalizing constant}}
\prod_{1\le i<j\le M}(\xi^M_i-\xi^M_j+j-i)^2\prod_{i=1}^M\f{(\xi^M_i+M)!}{(\xi^M_i+M-i)!}.\ee
It can easily be seen that the quantities of \eqre{2551419h} and \eqre{2551419h1} agree well at the limit.\epr

Let us now give a few definitions. Fix some positive integers $k, N$. 

\bgt\ite 
  An  \emph{up-right path  with values in  
  $\{1, \ld, N\}\ti\{1, \ld, k\}$} is a map $\pi$ defined on a set $\{1, \ld, p\}$ for a positive integer $p$ taking values in $\{1, \ld, N\}\ti\{1, \ld, k\}$, hence having two components $\pi_1, \pi_2$, \st  for all $i\in \{2, \ld, p\}$, we have the equality of sets $$\{\pi_1(i)-\pi_1(i-1),\pi_2(i)-\pi_2(i-1)\}=\{0,1\}.$$ 

\ite The  \emph{support} of an up-right path $\pi : \{1, \ld, p\}\to \{1, \ld, N\}\ti\{1, \ld, k\}$ is the set $\{\pi(1), \ld, \pi(p)\}$. 

\ite The  \emph{starting point} (resp. \emph{ending point}) of an up-right path $\pi : \{1, \ld, p\}\to \{1, \ld, N\}\ti\{1, \ld, k\}$ is     $\pi(1)$ (resp. $\pi(p)$). The  \emph{starting abscissa} (resp. \emph{ending abscissa}) of $\pi$ is the first coordinate of its  {starting point} (resp.  {ending point}).

\ite For $\pi$ an  up-right path  with values in $\{1, \ld, N\}\ti\{1, \ld, k\}$ and $P$ a point of the support of $\pi$, we denote by $\pi^{[P}$ (resp. $\pi^{P]}$) the path obtained by removing all points in $\pi$ strictly before (resp. strictly after) the point $P$ (note that there is no ampiguity in such a definition, as with our definition, an up-right path is always one-to-one).

\ite For $\pi, \pi^\prime$  up-right paths with values in $\{1, \ld, N\}\ti\{1, \ld, k\}$ \st the starting point of $\pi'$ is just on the right or just above the ending point of $\pi$, we denote by $\pi\cup \pi'$ their (obviously defined) concatenation.

\ite For $\pi, \pi^\prime$  up-right paths with values in $\{1, \ld, N\}\ti\{1, \ld, k\}$, we 
write $\pi<\pi'$, if for any $n\in \{1, \ld, N\}$, the intersection of the support 
of $\pi$ with $\{n\}\ti\{1, \ld, k\}$ is located strictly below the 
intersection of the support of $\pi^\prime$ with $\{n\}\ti\{1, \ld, k\}$.
\ent
We will also need the following lemmas.

\beg{lem}Let $1\le \ell \le k$ and let $\pi_1, \ld, \pi_\ell$ be some up-right paths with values in $\{1, \ld, N\}\ti\{1, \ld, k\}$ with pairwise disjoint supports. Then there exists   some up-right paths $\pi'_1, \ld, \pi'_{\ell}$ with values in $\{1, \ld, N\}\ti\{1, \ld, k\}$ with pairwise disjoint supports  and starting abscissas all equal to $1$ \st $$\cup_{i=1}^\ell \op{support}(\pi_i)\subset \cup_{i=1}^{\ell} \op{support}(\pi_i').$$\en{lem}

\bpr Note first that if $\ell=k$, the lemma is clearly true by choosing $\pi_1, \ld, \pi_k$ to be the maximal horizontal lines. So let us suppose that $\ell \le k-1$.  

Let us first state an obvious claim, which will be usefull in the sequel:
\\
{\bf Claim : } The hypothesis $\ell \le k-1$ implies that at least one point of the set $\{1\}\ti\{1, \ld, N\}$ is the starting point of no $\pi_i$ or that at least one $\pi_i$ has starting point in $\{1\}\ti\{1, \ld, N\}$ and second point just above its starting point.

Let $\op{sa}_{\max}$ denote the maximum of the starting abscissas of $\pi_1, \ld, \pi_\ell$ and $\ell_{\op{sa}_{\max}}$ denote the number of $i$'s \st the starting abscissa of $\pi_i$ is $\op{sa}_{\max}$. Let us prove the lemma by induction (for the lexical order) on $(\op{sa}_{\max},\ell_{\op{sa}_{\max}})$. If $\op{sa}_{\max}=1$, then there is nothing to do (as all $\pi_i$'s have starting abscissas equal to $1$). This starts the induction. 
So let us now suppose that $\op{sa}_{\max}\ge 2$. Let $i_0$ be \st $\pi_{i_0}$ has starting abscissa $\op{sa}_{\max}$. Let $P$ be the point in $\{1, \ld, N\}\ti\{1, \ld, k\}$ just on the left of the starting point of $\pi_{i_0}$. \bgt\ite If $P$ is on none of the supports of the $\pi_i$'s, then one can extend $\pi_{i_0}$ by adding $P$ at its beginning and conclude by the induction hypothesis.
\ite If $P$ is on the support of a path $\pi_{i_1}$, then two cases can occur:
\bgt\ite If $\pi_{i_1}$ does not end at $P$: then the point after $P$ in $\pi_{i_1}$ must be the point $Q$ just above $P$, and replacing $\pi_{i_0}$  by $\pi_{i_1}^{P]}\cup \pi_{i_0}$ and $\pi_{i_1}$ by $\pi_{i_1}^{[Q}$ allows to conclude by the induction hypothesis.
 \ite If $\pi_{i_1}$  ends at $P$: then replacing $\pi_{i_0}$ 
 by $\pi_{i_1}\cup \pi_{i_0}$ and $\pi_{i_1}$ by a point of the set $\{1\}\ti\{1, \ld, N\}$  allows to conclude  (note that replacing $\pi_{i_1}$ by a point of the set $\{1\}\ti\{1, \ld, N\}$  is possible by the Claim above : if point of the set $\{1\}\ti\{1, \ld, N\}$ is the starting point of no $\pi_i$, then this is obvious and if a $\pi_{i_2}$ has starting point in $\{1\}\ti\{1, \ld, N\}$ and second point just above its starting point, one can remove its  starting point to $\pi_{i_2}$ to get back to the first case).\ent
\ent
\epr

\beg{lem}Let $1\le \ell \le k$ and let $\pi_1, \ld, \pi_\ell$ be some up-right paths with values in $\{1, \ld, N\}\ti\{1, \ld, k\}$ with pairwise disjoint supports  and starting abscissas all equal to $1$. Then there exists   some up-right paths $\pi'_1, \ld, \pi'_{\ell}$ with values in $\{1, \ld, N\}\ti\{1, \ld, k\}$ with pairwise disjoint supports, starting abscissas all equal to $1$, ending abscissas all equal to $N$ \st $$\cup_{i=1}^\ell \op{support}(\pi_i)\subset \cup_{i=1}^{\ell} \op{support}(\pi_i').$$\en{lem}

\bpr  The proof is quite similar to the one of the previous lemma, but not exactly, as the fact that the starting abscissas are all equal to one are is now  added to the hypotheses and to the conclusion. 

Let $\op{ea}_{\min}$ denote the minimum of the ending abscissas of $\pi_1, \ld, \pi_\ell$ and $\ell_{\op{ea}_{\min}}$ denote the number of $i$'s \st the ending abscissa of $\pi_i$ is $\op{ea}_{\min}$. Let us prove the lemma by induction (for the lexical order) on $(\op{ea}_{\min},\ell_{\op{ea}_{\min}})$. If $\op{ea}_{\min}=N$, then there is nothing to do (as all $\pi_i$'s have ending  abscissas equal to $N$). This starts the induction. 
So let us now suppose that $\op{ea}_{\min}\le N-1$. Let $i_0$ be \st $\pi_{i_0}$ has ending abscissa $\op{ea}_{\min}$. Let $P$ be the point in $\{1, \ld, N\}\ti\{1, \ld, k\}$ just on the right of the ending point of $\pi_{i_0}$. \bgt\ite If $P$ is on none of the supports of the $\pi_i$'s, then one can extend $\pi_{i_0}$ by adding $P$ at its end and conclude by the induction hypothesis.
\ite If $P$ is on the support of a path $\pi_{i_1}$, then $P$ is not the starting point of $\pi_{i_1}$ and the point preceding  $P$ in $\pi_{i_1}$ is the point $Q$ just below $P$. Hence replacing $\pi_{i_0}$  by $ \pi_{i_0}\cup\pi_{i_1}^{[P}$ and $\pi_{i_1}$ by $\pi_{i_1}^{Q]}$ allows to conclude by the induction hypothesis.
\ent
\epr

\beg{lem}Let $1\le \ell \le k$ and let $\pi_1, \ld, \pi_\ell$ be some up-right paths 
with values in $\{1, \ld, N\}\ti\{1, \ld, k\}$ with pairwise disjoint supports, starting abscissas all equal to $1$ and  ending abscissas all equal to $N$. Then the  collection $(\pi_i)_{1\le i\le \ell} $ 
can be re-indexed in such a way that $\pi_1<\cd<\pi_\ell$.
\en{lem}

\bpr It suffices to re-index the  collection $(\pi_i)_{1\le i\le \ell} $ according to the orders of their starting abscissas and to use the definition of up-right paths and the  fact that their supports are disjoint. 
\epr

\beg{lem}\la{last}For each $1\le \ell\le k\le M$, \be\la{micheltalomorvan}
\lam_1^k+\cd+\lam_\ell^k=\max_{\pi_1, \ld,\pi_\ell}\sum_{r=1}^\ell \sum_{(i,j)\in \pi_r} w_{ij} ,\ee
where the max is over collections $\{\pi_1<\cd<\pi_\ell\}$ of  up-right paths 
in the set $$\{1, \ld, N\}\ti\{1, \ld, k\}$$  starting in the 
subset  $\{1\}\ti \{1, \ld, k\}$  and ending in the subset $\{N\}\ti \{1, \ld, k\}$. \en{lem}

\bpr By   \cite[Th.~1.1.1]{moerbeke} and \cite[Chap.~3, Lemma 1]{fulton},
  \eqre{micheltalomorvan} is true when the paths are only required to be 
pairwise disjoint, without any condition on the starting and ending points.  Then, the three  previous lemmas allow to claim that any set of $\ell$ 
pairwise disjoint paths can be changed into a set of $\ell$ pairwise disjoint paths $\pi_1<\cd<\pi_\ell$, with starting abscissas all equal to $1$ and  ending abscissas all equal to $N$
  in such a way that the union of the  supports of the news paths contains the union of the supports of the former ones.      To finish the proof, it suffices then 
to notice that since the $w_{ij}$'s are non-negative, enlarging the union of the supports of the 
paths never decreases the total weight.
\epr

To complete the proof of Theorem~\re{208121}, note that any up-right path $\pi_r$ as described 
in the previous lemma is a 
concatenation of at most $k$ paths with fixed second coordinate and has length 
between $N$ and $N+M$.  Moreover, by Donsker theorem (see \cite{billingsley,GW}), the 
$M$-dimensional process
$$\lf(\ff{\sqrt{vN}}\sum_{i=1}^{\lfloor Nt\rfloor}(w_{ik}-e)\ri)_{k=1, \ld, M} $$ converges in 
distribution (for the Skorohod topology) to the 
$M$--dimensional Brownian motion $B$.  To finish the proof, apply both Lemma~\re{1481218h07ter}  
and Lemma~\re{last} .
\epr

\section{Applications of Theorem~\re{208121}}
 
Let us present and prove at first some of its corollaries which  
motivated, in large part, the present study.   When combined with \cite{HL1}, the first corollary 
provides an alternative approach to \cite{TW} or \cite{ITW1, ITW2}.  
The second corollary makes full use of Theorem~\re{208121} and, when combined 
with \cite{HX} or \cite{HL2}, provides an alternative 
approach to \cite{KJAOM2001}, \cite{ITW1}, \cite{ITW2} (see also Bufetov\cite{Bu} and M\'eliot\cite{Me}).

Let us briefly recall the framework of the works just cited.  
Let $(X_i)_{i\ge 1}$ be a sequence 
of i.i.d.~random variables on a totally ordered finite alphabet $\A$ of cardinality $k$.  
Denote the elements of $\A$ by $\al_1, \ld, \al_k$ listed in such a way  
that if $p_i:=\p(X_1=\al_i)$, $i=1, \ld, k$, then $p_1\ge \cdd\ge p_k$ 
(therefore this indexing of the letters in $\A$ has nothing to do 
with the order used on $\A$).  Next, decompose the alphabet 
$\A$ into subsets $\A_1, \ld, \A_n$  in such a way 
that $\al_i$ and $\al_j$ belong to the same $\A_m$, $m=1,\ld, n\le k,$ \ssi $p_i=p_j$.  
Finally, let $\op{LI}_N$ be the length of the longest increasing 
subsequence of the random word $$X_1\cd\cd X_N.$$

\beg{cor} Let $p_{\max}:=p_1$, $k_1:=\# \A_1$ and let $H=[h_{ij}]$ be a 
$k_1 \times k_1$ GUE matrix 
with largest eigenvalue $\mu_{\max}$.  Then, as 
$N$ tends to infinity, 
\begin{align}\label{convbm1} 
&\frac{\op{LI}_N-Np_{\max}}{\sqrt{Np_{\max}}} \Longrightarrow
     \frac{\sqrt{1-k_1p_{\max}}-1}{k_1}
     \sum^{k_1}_{j=1} h_{jj} + \mu_{\max}.
\end{align}
 \en{cor}

\bpr 
From Theorem~\re{208121}, with the notation introduced above, and if 
$B$ is now a $k_1$--dimensional standard Brownian motion, 
\begin{equation}\la{eqlaw1}
(\mu_1^j)_{1\le j\le k_1}\cup(h_{jj})_{1\le j\le k_1}\;\;\!\eqlaw \;\;\!\left(\mathop{\max_{
    {\scriptstyle 0 = t_0 \le  \cdots\le t_j = 1}}}  
    \sum_{i=1}^j(B_i({t_i})-B_i({t_{i-1}}))\!\right)_{1\le j\le k_1}\!\!\!\!\!\!\!\cup\;\;(B_j(1))_{1\le j\le k_1},  
\end{equation}
since $h_{11}=\mu_1^1$ and for all $j=2, \ld, k_1$,  
$h_{jj} = \sum_{i=1}^j \mu_{i}^j - \sum_{i=1}^{j-1}\mu_i^{j-1}$.  
Next, Corollary 3.3 in \cite{HL1} asserts that 
\begin{equation}\label{convbm2}
\frac{\op{LI}_N-Np_{\max}}{\sqrt{Np_{\max}}}  \Longrightarrow
     \frac{\sqrt{1-k_1p_{\max}}-1}{k_1}
     \sum^{k_1}_{j=1} B_{j}(1) + \mathop{\max_{
    {\scriptstyle 0 = t_0 \le  \cdots\le t_{k_1} = 1}}}
    \sum^{k_1}_{i=1}  \left( B_{i}(t_{i})
          - B_{i}(t_{i-1}) \right),
\end{equation}
combining \eqre{eqlaw1} and \eqre{convbm2} gives \eqre{convbm1}.
\epr

Denote by $\lam_1\ge \cd\ge \lam_k$ the shape of the Young diagrams 
obtained by applying the RSK algorithm to the random 
word $$X_1\cd\cd X_N,$$ 
and let  
$$\xi_i=\f{\lam_i-Np_i}{\sqrt{Np_i}},$$
$1 \le i \le k$, be the corresponding rescaled variables.  
Introduce next some independent GUE 
matrices $H_1, \ld, H_n$, where each $H_j$ has 
size $k_j:=\# \A_j$, and let $$H:=\bpm H_1&&\\
&\ddots&\\
&&H_n\epm\qquad\trm{ and }\qquad \bpm \wtH_1&&\\
&\ddots&\\
&&\wtH_n\epm:=H-\Tr(HJ)J,$$ where $$J=\diag(\sqrt{p_1}, \ld, \sqrt{p_k}).$$
\beg{rmk}{\rm Note that $J$ is a unit vector of the 
space $\h_{k_1}\ti \cd\ti \h_{k_n}$ endowed with 
the Euclidean product structure, so $H-\Tr(HJ)J$ is the orthogonal projection  
onto $J^\perp$, so that its law is the law of $H$ conditioned to 
belong to $J^\perp$.}\en{rmk}

Finally, define the random vector $(\mu_1, \ld, \mu_k)$ 
by $$(\mu_1, \ld, \mu_k):=(\trm{ordered spectrum of }\wtH_1, \ld\ld, 
\trm{ordered spectrum of } \wtH_n).$$
 
\beg{cor}As $N\to\infty$, $$(\xi_1, \ld, \xi_k) \Longrightarrow (\mu_1, \ld, \mu_k). $$
\en{cor}

\beg{rmk}{\rm The limiting law of $\op{LI}_N$, rescaled, is simply the 
law of $\mu_1$ and is given by 
\begin{equation}\label{anotherrep}
\lam_{\max}(H_1)-p_1\Tr(H_1)-\sqrt{p_1(1-k_1p_1)}Z,
\end{equation}
where $Z$ is a standard Gaussian random variable, independent of $H_1$.  
Note also that this law only depends 
on $p_1$ and $k_1$.}\en{rmk}

\bpr  First, $$\Tr(HJ)=\sum_{j=1}^n \sqrt{p_{(j)}}\Tr H_j,$$ where for 
all $j$, $p_{(j)}:=p_\ell$ for $\ell\in \A_j$. 
So for each $i$, we have $$\wtH_i=H_i-\lf(\sqrt{p_i}\sum_{j=1}^n \sqrt{p_{(j)}}\Tr H_j\ri)I,$$  
where $I$ is the corresponding identity matrix.  
Then, Theorem~3.1 and Remark~3.2 (iv) in~\cite{HX} 
together with Theorem~\re{208121} allow to conclude.\epr

In case the i.i.d.~random variables generating the random word are replaced 
by an (irreducible, aperiodic) 
homogeneous Markov chain, with state space $\A$ of cardinality $k$, 
the corresponding limiting laws are also given in terms 
of maximal Brownian functionals similar to those 
in Theorem \re{208121} (see \cite{HL2}).  
However, an important difference is that now the standard Brownian motion $B$  
is replaced by a correlated one $\tilde B$ with, say, covariance matrix $\Sigma$ 
instead of $I$.   
The possible identification of (the law of) these functionals as (the law of) maximal 
eigenvalues (or spectra) of random matrices has not been fully accomplished yet, although various 
cases are done.  
In particular, for cyclic transition matrices $P$, in which case the stationary distribution is 
the uniform one, there is a curious dichotomy 
between alphabets of size at most three and size four or more.  
Indeed for $k\le 3$, the cyclic hypothesis forces $\Sigma$ to have a 
permutation-symmetric structure seen in the i.i.d.~uniform case.  For example, 
for $k=3$, $\Sigma$ is, a rescaled version of, 
\begin{equation}\label{sigma4u111213}
{\Sigma_u}:= \begin{pmatrix}
     1      &-1/2   &-1/2\\
    -1/2   &1       &-1/2\\
    -1/2   &-1/2    &1
  \end{pmatrix},
\end{equation}
and so (up to a multiplicative constant) and with $k=k_1=3$, $p_{\max} = p_1 = 1/3$, 
\eqref{convbm1} continues to hold for cyclic Markov chains.   
For $k\ge 4$, the cyclicity constraint on $P$ forces $\Sigma$ to be cyclic but does no 
longer force the permutation-symmetric structure, and, say, for $k=4$, $\Sigma$ might differ 
from, a rescaled version of, 
\begin{equation}\label{sigma4u}
{\Sigma_u}:=\begin{pmatrix}
    1      &-1/3   &-1/3   &-1/3\\
    -1/3   &1      &-1/3   &-1/3\\
    -1/3   &-1/3   &1      &-1/3\\
    -1/3   &-1/3   &-1/3   &1
  \end{pmatrix}.  
\end{equation}         
In fact, if 
\begin{equation}\label{matrixp4}
P =
  \begin{pmatrix} 
  p_1 & p_2 & p_3 & p_4\\
  p_4 & p_1 & p_2 & p_3\\
  p_3 & p_4 & p_1 &p_2\\
  p_2 & p_3 & p_4 & p_1
  \end{pmatrix}, 
\end{equation}
then $\Sigma$ is a rescaled version of $\Sigma_u$ if and only if $p_3^2=p_2p_4$.  
Nevertheless, see \cite{HL2},   for $k\ge 2$,  
\begin{align}\label{conv}
\frac{\op{LI}_N -N/k}{\sigma\sqrt{N}}  \Longrightarrow  
\sup\!\left\{\!\Del_\pi(\tilde B); \pi\in\Pa , \pi\le k\!\right\}
=  \max_{0=t_0\le  \cdots\le t_k=1} 
      \sum_{j=1}^{k}\!\left(\!\tilde B_{j}(t_{j})  - \tilde B_{j}(t_{j-1})\!\right).   
\end{align}

Assuming that in addition to be cyclic, $P$ is also symmetric (for $k=2$ the cyclic and symmetric 
assumptions are the same, and 
Chistyakov and G\"otze \cite{CG}, see also \cite{HL3}, showed that the 
corresponding limiting law is the maximal eigenvalue of the $2\times 2$ traceless GUE) 
a diagonalization argument, combined with \eqref{conv}, leads to the 
following result. 

\beg{propo}
Let $P:= (p_{i,j})_{1\le i,j\le k}$ be cyclic and symmetric, i.e., 
$P=  (p(j-i))_{1\le i,j\le k}$, where $p$ is a 
$k$-periodic function defined on $\z$ such that $p(r) = p(-r)$, for all $r\in\z$.  
Let \be\la{112131}\lambda_\ell := \sum_{r=1}^kp(r)\cos\left({2\pi(\ell-1)r}/{k}\right)\,\qquad \trm{ $(1\le \ell \le k)$},\ee 
and let $(B_j)_{j=2,\ldots,k}$ be a $(k-1)$-dimensional standard 
Brownian motion on $[0,1]$.  Then, 
\begin{eqnarray}\nonumber 
\frac{\op{LI}_N -N/k}{\sigma\sqrt{N}} \Longrightarrow \max_{0=t_0\le  \cdots\le t_k=1}  
\left\{\!\!\sqrt{\frac{2}{k}}\!\sum_{j=1}^k\!\!\sum_{r=1}^{\left\lfl\frac{k-1}{2}\right\rfl}
\!\sqrt{\frac{1+\lambda_{r+1}}{1-\lambda_{r+1}}}
\!\left(\!\cos\!\left(\!\frac{2\pi jr}{k}\!\right)\!(B_{2r}(t_j) - \right.\right.\\
\quad \quad \quad \quad \left. B_{2r}(t_{j-1})) + \sin\!\left(\!\frac{2\pi jr}{k}\!\right)(B_{2r+1}(t_j) 
- B_{2r+1}(t_{j-1}))\!\right) \nonumber \\ \label{firstprop}
\quad \quad \left.-{\frac{1}{\sqrt k}}
\sqrt{\frac{1+\lambda_{\frac{k}{2}+1}}
{1-\lambda_{\frac{k}{2}+1}}} \left(2\sum_{j=1}^k \left(B_k(t_j)-B_k(t_{j-1})\right)-B_k(1)\right)\right\}, 
\end{eqnarray}
where the last term above is only present for $k$ even.  
\en{propo} 
 
\bpr
Since $P$ is symmetric, it can be diagonalized as 
$P=S\Lambda S\tra$, where $\Lambda$ is the diagonal matrix formed with its 
eigenvalues $(\lambda_\ell)_{1\le \ell \le k}$ (we will see below that these are the 
quantities defined at \eqre{112131}) and where $S$ is a matrix formed by the orthonormal 
column eigenvectors $(u_\ell)_{1\le\ell\le k}$ where 
$u_1\tra = (1/\sqrt k, \dots, 1/\sqrt k)$.  Next, by Theorem~4.3 in \cite{HL2}, $\Sigma$, the 
covariance matrix of the $k$--dimensional correlated Brownian motion $\tilde B$, is given by  
$\Sigma = S {\Lambda_{\Sigma}}S\tra$, where $\Lambda_{\Sigma}$ is the diagonal matrix 
with diagonal entries $0, {(1+\lambda_2)}/{(1-\lambda_2)},\ld, {(1+\lambda_k)}/{(1-\lambda_k)}$.  
Therefore, $\tilde B = S\sqrt{\Lambda_\Sigma} B$, where now $B$ is a 
standard $k$--dimensional Brownian motion.   
Now the symmetric and cyclic structures imply that the eigenvalues of $P$ are 
in fact the $\lam_\ell$'s defined at \eqre{112131}:  
$\lambda_\ell = \sum_{r=1}^kp(r)\cos\left({2\pi(\ell-1)r}/{k}\right)$, $1\le \ell \le k$, 
(clearly they are not all simple since  
$\lambda_\ell = \lambda_{k-\ell+2}, \ell=2, \dots, k$).     
The corresponding orthonormal column eigenvectors are 
$$v_\ell :=(v_{j,\ell})_{1\le j \le k}=(\sqrt{2}\cos(2\pi(\ell-1)j/k)/\sqrt k)_{1\le j \le k}\,,\qquad 
\ell = 1, 2, \dots, \lfl k/2\rfl+1,$$ and 
$$w_\ell :=(w_{j,\ell})_{1\le j \le k}
=(\sqrt 2\sin(2\pi(\ell-1)j/k)/\sqrt k)_{1\le j \le k}\,,\qquad 
\ell = 2, 3, \dots, \lfl(k-1)/2\rfl+1.$$           
Clearly, $v_1=u_1$ is an eigenvector corresponding to the simple eigenvalue $1$, while if $k$ is 
even, $v_{(k/2)+1} = (1/\sqrt k, -1/\sqrt k, \dots, 1/\sqrt k,-1/\sqrt k)$ is an eigenvector 
corresponding to the simple eigenvalue $\sum_{r=1}^k p(r)\cos\left(2\pi(k/2+1-1)r/k\right) = \sum_{r=1}^k(-1)^rp(r)$.  Moreover, 
for $\ell = 2, 3, \dots, \lfl(k-1)/2\rfl+1$, $v_\ell$ and $w_\ell$ share 
the same eigenvalue $\lambda_\ell$.  
Therefore,  
$$
S= \left(v_1, v_2, w_2, \dots, v_{\lfl\frac{k-1}{2}\rfl+1}, w_{\lfl\frac{k-1}{2}\rfl+1}, v_{\frac{k}{2}+1} \right), 
$$  
where, above, the last column is only present if $k$ is even.  
Next, from the transformation $\tilde B = S\sqrt{\Lambda_\Sigma} B$, and since 
$\left(\sqrt{\Lambda_{\Sigma}} B\right)_\ell = \sqrt{(1+\lambda_{\lfl\ell/2\rfl+1})/(1-\lambda_{\lfl\ell/2\rfl+1})}B_\ell$, 
$\ell=2, \dots, k$, and $\left(\sqrt{\Lambda_{\Sigma}} B\right)_1=0$, then for $j=1, \dots, k$, 
\begin{align}\label{conv2}
\tilde B_j = \sum_{\ell=2}^{k} u_{j,\ell}\sqrt{\frac{1+\lambda_{\left\lfl\frac{\ell}{2}\right\rfl+1}}
{1-\lambda_{\lfl\frac{\ell}{2}\rfl+1}}}B_\ell,   
\end{align}  
where $u_{j,\ell} = v_{j,\lfl\ell/2\rfl+1}$ or $u_{j,\ell}=w_{j,\lfl\ell/2\rfl+1}$, for $\ell$ even or odd. 
Therefore, for $j = 1, \dots, k$, 
\begin{align}\label{conv3}
\tilde B_j = \sqrt{\frac{2}{k}}\sum_{r=1}^{\left\lfl\frac{k-1}{2}\right\rfl} 
\sqrt{\frac{1+\lambda_{r+1}}
{1-\lambda_{r+1}}} &\left( \cos\left(\frac{2\pi rj}{k}\right)B_{2r} + 
\sin\left(\frac{2\pi rj}{k}\right)B_{2r+1}\right) \\
& \quad \quad + {\frac{(-1)^{j+1}}{\sqrt k}}
\sqrt{\frac{1+\lambda_{\frac{k}{2}+1}}
{1-\lambda_{\frac{k}{2}+1}}} B_k, \nonumber   
\end{align}  
where the last term on the right of \eqref{conv3} is only present for $k$ even.  
With \eqref{conv3}, the sum on the right hand 
side of \eqref{conv} becomes:  
\begin{align}\label{conv4}
\sqrt{\frac{2}{k}}\sum_{j=1}^k\sum_{r=1}^{\left\lfl\frac{k-1}{2}\right\rfl}
\sqrt{\frac{1+\lambda_{r+1}}{1-\lambda_{r+1}}}
\left(   \cos\left(\frac{2\pi jr}{k}\right) B_{2r}(t_j) +  
\sin\left(\frac{2\pi jr}{k}\right)B_{2r+1}(t_j)\right.\\ 
\quad \quad  \quad \left. - \cos\left(\frac{2\pi jr}{k}\right) B_{2r}(t_{j-1}) 
- \sin\left(\frac{2\pi jr}{k}\right) B_{2r+1}(t_{j-1})\right) \nonumber \\
+{\frac{1}{\sqrt k}}
\sqrt{\frac{1+\lambda_{\frac{k}{2}+1}}
{1-\lambda_{\frac{k}{2}+1}}} \sum_{j=1}^k (-1)^{j+1}\left(B_k(t_j)-B_k(t_{j-1})\right), 
\nonumber
\end{align}
an expression only involving standard Brownian motions and where, again, the last term 
\begin{align}\label{conv5}
{\frac{1}{\sqrt k}}
\sqrt{\frac{1+\lambda_{\frac{k}{2}+1}}
{1-\lambda_{\frac{k}{2}+1}}}&\sum_{j=1}^k (-1)^{j+1}\left(B_k(t_j)-B_k(t_{j-1})\right)\\
&\quad \quad = {\frac{-1}{\sqrt k}}
\sqrt{\frac{1+\lambda_{\frac{k}{2}+1}}
{1-\lambda_{\frac{k}{2}+1}}} \left(\!\!B_k(1) + \sum_{j=1}^{k-1} 2(-1)^{j}B_k(t_j)\!\!\right)\nonumber\\
&\quad \quad = {\frac{-1}{\sqrt k}}
\sqrt{\frac{1+\lambda_{\frac{k}{2}+1}}
{1-\lambda_{\frac{k}{2}+1}}} \left(2\sum_{j=1}^{k}(B_k(t_j) - B_k(t_{j-1}))- B_k(1)\right),  
\end{align} 
is only present if $k$ is even.

\epr
 
\beg{rmk}{\rm
Let us try to specialize the previous results in instances where further simplifications 
and identifications occur.  

(i) For $k=3$, and up to the factor 
$\sqrt{2(1+\lambda_{2})/(k(1-\lambda_{2}))}=\sqrt{{2(1+3p_1)}/{(3(3-3p_1))}}$, the right-hand side of \eqref{firstprop} becomes 
\begin{align*}
\max_{0=t_0\le t_1 \le t_{2}\le t_3=1}& 
\sum_{j=1}^3 \!\left(\!\cos\!\left(\!\frac{2\pi j}{3}\!\right)\!(B_{2}(t_j) -  B_{2}(t_{j-1})) + 
\sin\!\left(\!\frac{2\pi j}{3}\!\right)(B_{3}(t_j) - B_{3}(t_{j-1}))\!\right) \\
&=\max_{0\le t_1 \le t_{2}\le 1}\left(B_2(1) + \sqrt{3}B_3(t_1)-\frac{\sqrt{3}}{2}
B_3(t_2)-\frac{3}{2}B_2(t_2)\right)\\
&\eqlaw \sqrt{\frac{2}{3}}\max_{0\le t_1 \le t_{2}\le t_3=1} 
\sum_{j=1}^2 \!\left(-\sqrt{\frac{j}{j+1}} B_{j}(t_{j+1}) + \sqrt{\frac{j}{j+1}} B_{j}(t_{j}) \right),\\ 
&\eqlaw \sqrt{\frac{2}{3}}\left(\max_{0=t_0\le t_1 \le t_{2}\le t_3=1} 
\sum_{j=1}^3 \!\left( B_{j}(t_j) -  B_{j}(t_{j-1})\right) -\frac{1}{3}\sum_{j=1}^{3}B_j(1)\right),  
\end{align*}
where the last equality, in law, follows 
either by using, in \eqref{conv}, 
the simple linear transformation 
\begin{equation}\label{trans1}
\tilde B_j = 
\sqrt{\frac{2(1+\lambda_{2})}{3(1-\lambda_{2})}}\left(\sqrt{\frac{2}{3}}B_j - \sqrt{\frac{1}{6}}\sum_{i=1, i\neq j}^3B_i\right), \quad j =1,2,3,
\end{equation} 
which, by comparing covariances, is easily verified from 
\eqref{conv3}; or, still by comparing covariances, by arguments such as those 
in the proof of Theorem~3.2 in \cite{HL1}.  
Therefore, with the help of Theorem~\re{208121}, and up to a scaling factor, 
the limiting law of $\op{LI}_N$ is that of the maximal eigenvalue of the 
$3\ti 3$ traceless GUE.

(ii) For $k=4$, $\Sigma$, the covariance 
matrix of $\tilde B = (\tilde B_j)_{j=1, \dots, 4}$ is given, up to a scaling constant, by: 
\begin{equation}\label{sigma4}
{\Sigma}:=\begin{pmatrix}
   2\eta_2 + \eta_3       &-\eta_3   &-2\eta_2+\eta_3   &-\eta_3\\
    -\eta_3   &2\eta_2 + \eta_3      &-\eta_3   &-2\eta_2+\eta_3\\
    -2\eta_2+\eta_3  &-\eta_3   &2\eta_2 + \eta_3      &-\eta_3\\
      -\eta_3 &-2\eta_2+\eta_3   &-\eta_3   &2\eta_2 + \eta_3
  \end{pmatrix}. 
\end{equation}
where $\eta_2 = (1+\lambda_2)/(1-\lambda_2)$, $\lambda_2 = p_1 - p_3$, 
and $\eta_3 = (1+\lambda_3)/(1-\lambda_3)$, $\lambda_3 = p_1 -2p_2 + p_3$.  
Clearly, $\Sigma$ can differ from $\Sigma_u$, e.g., 
let $2\eta_2 = \eta_3$, i.e., 
let    

\begin{equation*}
P =
  \begin{pmatrix} 
  p_1 & p_2 & \frac{p_2(1-2p_2)}{1+2p_2} & p_2\\
  p_2 & p_1 & p_2 & \frac{p_2(1-2p_2)}{1+2p_2}\\
  \frac{p_2(1-2p_2)}{1+2p_2} & p_2 & p_1 &p_2\\
  p_2 & \frac{p_2(1-2p_2)}{1+2p_2} & p_2 & p_1
  \end{pmatrix}.  
\end{equation*}  
Then, and up to the multiplicative constant $4\eta_2$, $\Sigma$ becomes:  
\begin{equation*}
{\Sigma}:=\begin{pmatrix}
    1       &-1/2   &0   &-1/2\\
    -1/2   &1      &-1/2   &0\\
    0  &-1/2   &1      &-1/2\\
      -1/2 &0   &-1/2   &1
  \end{pmatrix}, 
\end{equation*}
which is clearly different from, a rescaled version of, \eqref{sigma4u}.  
In fact, if $\Sigma = \Sigma_u$, then clearly   
\begin{equation}\label{trans2}
\tilde B_j = 
\frac{\sqrt 3}{2}B_j - \frac{1}{2\sqrt3}\sum_{i=1, i\neq j}^4B_i, \quad j =1,2,3,4.  
\end{equation}
Conversely, and as easily seen, for a linear transformation such as  
\begin{equation*}
\tilde B_j = 
\alpha_jB_j - \sum_{i=1, i\neq j}^4\beta_iB_i, \quad j =1,2,3,4, 
\end{equation*} to lead to $\Sigma$, one needs to $\Sigma$ to be 
  permutation-symmetric  and, up to a multiplicative constant, 
the right-hand side of \eqref{firstprop} becomes equal in law to  
$$\max_{0=t_0\le t_1 \le t_{2}\le t_3\le t_4=1} 
\sum_{j=1}^4 \!\left( B_{j}(t_j) -  B_{j}(t_{j-1})\right) -\frac{1}{4}\sum_{j=1}^{4}B_j(1),   
$$
and corresponds to the matrix $P$ in \eqref{matrixp4} with $p_2=p_3=p_4$.  

(iii) Finally, it is easy to see that the properties just described 
continue to hold for arbitrary dimension $k\ge 4$.  In arbitrary dimension, 
if $\Sigma = \Sigma_u$ (the $k$-dimensional version of the matrix defined at \eqre{sigma4u111213} and \eqre{sigma4u}), then the    
linear transformation corresponding to \eqref{trans2} is given by   
\begin{equation*}
\tilde B_j = 
\sqrt{\frac{k-1}{k}}B_j - \sqrt{\frac{1}{k(k-1)}}\sum_{i=1, i\neq j}^kB_i, \quad j =1,\dots,k, 
\end{equation*} 
Conversely,  for a linear transformation such as  
\begin{equation*}
\tilde B_j = 
\alpha_jB_j - \sum_{i=1, i\neq j}^k\beta_iB_i, \quad j =1,2,\dots,k,  
\end{equation*} to lead to $\Sigma$, one needs to $\Sigma$ to be 
  permutation-symmetric.  
In either instance, and up to a multiplicative constant, 
the right-hand side of \eqref{firstprop} has the same law as  
$$\max_{0=t_0\le   \cdots \le t_k=1}  
\sum_{j=1}^k \!\left( B_{j}(t_j) -  B_{j}(t_{j-1})\right) -\frac{1}{k}\sum_{j=1}^{k}B_j(1),   
$$
which, in turn, via Theorem~\re{208121}, is equal in law to the maximal eigenvalue of an 
element of the $k\ti k$ traceless GUE.}   
\end{rmk}

\begin{thebibliography}{10}
 \bibitem{alice-greg-ofer}
G.~Anderson, A.~Guionnet, O.~Zeitouni
  \emph{An Introduction to Random Matrices}.
  Cambridge studies in advanced mathematics, {118} (2009).

\bibitem{baryshnikov} Y. Baryshnikov \emph{GUEs and queues}, Probab. Theory Relat. Fields 119, 256--274 (2001).

\bibitem{BBO} P. Biane, P. Bougerol, N. O'Connell \emph{Littelmann paths and Brownian paths},
{Duke Math. J.}, vol. 130, pp. 127--167, 2005.

\bibitem{billingsley} P. Billingsley \emph{Convergence of probability measures} Wiley, 1999.

\bibitem{BJ} P. Bougerol, T. Jeulin \emph{Paths in Weyl chambers and random matrices}
{Probab. Theory Related Fields}, vol. 124, no.~4, pp. 517--543, 2002.

\bibitem{Bu} A. Bufetov \emph{A central limit theorem for the 
extremal characters of the infinite symmetric group}, {Functional Analysis and Its Applications}, vol.~46, no.~2, 
pp.~83-93, 2012.

\bibitem{CG} G.P.~Chistyakov, F.~G\"otze \emph{Distribution of the shape of Markovian random words}, 
{Probab. Theory Relat. Fields}, vol. 129, no.~1, pp.~18-36, 2004.

\bibitem{doumerc} Y. Doumerc \emph{A note on representations of classical Gaussian matrices} { S\'emimaire de Probabilit\'es XXXVII.}, { Lecture Notes in Math., No. 1832}, Springer, Berlin, (2003),  370--384.

\bibitem{fulton} W. Fulton \emph{Young Tableaux}, CUP, Cambridge, (1997).

\bibitem{GW} P.~Glynn, W.~Whitt \emph{Departures form many queues in series}, Ann. Appl. 
Probability, 1, 546--572, 1991.

\bibitem{GTW} J.~Gravner, C.A.~Tracy, H.~Widom \emph{Limit theorems for height fluctuations in a class 
of discrete space and time growth models},  J.~Stat.~Phys.~102, (2001), 1085--1132.

\bibitem{HL1} C. Houdr\'e, T. Litherland \emph{On the longest increasing subsequence for finite and 
countable alphabets}, in {High Dimensional Probability V: The Luminy Volume} 
(Beachwood, Ohio, USA: Institute of Mathematical Statistics), 185--212, 2009.

\bibitem{HL2} C. Houdr\'e, T. Litherland \emph{On the limiting shape of Young diagrams associated 
with Markov random words}, arXiv 1110.4570, (2011).

\bibitem{HL3} C. Houdr\'e, T. Litherland \emph{Asymptotics for the length of the longest 
increasing subsequence of a binary Markov random word}, in: {Malliavin calculus and 
stochastic analysis: A Festschrift in honor of David Nualart}, New York, NY: Springer. 
Springer Proceedings in Mathematics \& Statistics 34, 511-524, 2013.

\bibitem{HX} C. Houdr\'e, H. Xu \emph{On the limiting shape of Young diagrams associated with 
inhomogeneous random words}, in: {High Dimensional Probability VI: The Banff volume} 
Progress in Probability, 66, Birkhauser (2013)

\bibitem{ITW1} A.~Its, C.A.~Tracy, H.~~Widom \emph{Random words, Toeplitz determinants, and 
integrable systems. I.} Random matrix models and their 
applications, 245--258, Math. Sci. Res. Inst. Publ., 40, Cambridge Univ. Press, Cambridge, 2001.

\bibitem{ITW2} A.~Its, C.A.~Tracy, H.~Widom \emph{Random words, Toeplitz determinants, and 
integrable systems. II.}  
Advances in nonlinear mathematics and science. Phys. D., vol. 152-153 (2001), 199--224.

\bibitem{J1} K. Johansson \emph{Shape fluctuations and random matrices} 
Comm. Math. Phys. 209, (2000), 437--476.

\bibitem{KJAOM2001} K. Johansson \emph{Discrete orthogonal polynomial ensembles 
and the Plancherel measure}.
Ann. of Math. (2) 153 (2001), no. 1, 259--296.

\bibitem{J2} K. Johansson \emph{Non-intersecting paths, random tilings and 
random matrices}, Probab.~Theor.~Rel.~Fields, 123, (2002), 225--280.

\bibitem{NJEJP06} K. Johansson, E. Nordenstam \emph{Eigenvalues of GUE minors}, 
Elec. J. Probab. Vol. 11 (2006), Paper no. 50, pages 1342--1371.

\bibitem{kuperberg2002} G. Kuperberg \emph{Random words, quantum statistics, 
central limits, random matrices}. Methods Appl. Anal. 9 (1), 99--118 (2002).  

\bibitem{Me} P.-L. M\'eliot \emph{Fluctuations of central measures on partitions} 
{Proceedings of the 24th International Conference on Formal Power series and Algebraic Combinatorics} 
(Nagoya, Japan), pp.~387-398, 2012.

\bibitem{moerbeke} P. van Moerbeke \emph{Random and Integrable Models in Mathematics
and Physics} in  \emph{Random Matrices, Random Processes and Integrable Systems}, 
edited by J. Harnad, CRM Series in mathematical Physiscs,  2011.

\bibitem{O} N. O'Connell \emph{A path transformation for random walks and the 
Robinson-Schensted correspondence} Trans. Amer. Math. Soc. 355, (2003), 3669--3697.

\bibitem{OY} N. O'Connell, M.~Yor \emph{A representation for non-colliding random walks}
{Electron. Comm. Probab.}, vol. 7,  (2002), 1--12.

\bibitem{TW} C. Tracy, H. Widom \emph{On the distribution of the lengths of the longest increasing monotone  subsequences in random words}. {Probab. Theor. Rel. Fields.} {119} (2001), 350--380.

 \en{thebibliography}

\end{document}